\documentclass[]{{birkart}}

\usepackage{graphicx}
\usepackage{amsfonts}
\usepackage{amsmath,amscd}
\usepackage{amsthm}

{ \theoremstyle{plain}
\newtheorem{theorem}{Theorem}
\newtheorem{lemma}{Lemma}
\newtheorem{corollary}{Corollary}
\newtheorem{proposition}{Proposition}
}

{ \theoremstyle{definition}
}

{ \theoremstyle{definition}

}

{ \theoremstyle{definition}

}

\def\scal#1#2{\langle #1, #2\rangle}
\def\M{\mathcal{M}}

\def\R#1{\mathbb{R}^{#1}}

\def\D{{\overline{\nabla}}}

\DeclareMathOperator{\conv}{\mathrm{conv}}

\def\overt{\overline{t}}

\begin{document}


\title[External geometry of $p$-minimal surfaces]{External geometry of $p$-minimal surfaces}
\author{Tkachev Vladimir G.
}
\thanks{The paper was supported by RFRF,
project 93-011-176}

\begin{abstract}
A surface $\M$ is called $p$-minimal if one of the coordinate
functions is $p$-harmonic in the inner metric. We show that in the
twodimensional case the Gaussian map of such surfaces is
quasiconformal.  In the case when the surface is a tube we study the
geometrical structure of such surfaces. In particularly, we establish
the second order differential inequality for the form of the
sections of $\M$ which generalizes the known ones in the minimal
surfaces theory.
\end{abstract}

\maketitle

\section{Introduction}
\label{sec1}

{\bf 1.1.}
Let $\M=(M;x)$ be a surface given by a $C^2$-immersion $x: M \to
\R{n+1}$ of $n$-dimensional orientable noncompact manifold $M$.

\smallskip
{\bf Definition 1.}
A surface $\M$ is said to be {\it minimal} if its mean curvature
vector $H(m)\equiv 0$.

The well-known property of the minimal immersions in the Euclidean space
is harmonicity of their coordinate functions. Moreover, if
{\it one} coordinate function of an immersion is harmonic then {\it
all} coordinates satisfy this property and the immersion is minimal.
On the other hand, for $n=2$ this condition yields the fact that the
Gauss map of such surfaces is conformal \cite{Os}.

The natural question arises: what happens if we replace the
requirement of harmonicity by $p$-harmonicity?

\smallskip
{\bf Definition 2.} For a fixed $p>1$ a surface $\M$ is said to be
$p$-{\it minimal} if one of the coordinate functions is $p$-harmonic with
respect to the inner metric of $\M$. Other words, there exists a
direction $e\in \R{n}$ such that
\begin{equation}
\Delta_p f\equiv{\rm div}|\nabla f|^{p-2}\nabla f=0,
\label{plaplace}
\end{equation}
where $f(m)=\scal{x(m)}{e}$ and $\nabla $ is the covariant derivative
on $M$.

One easily can be shown that $p$-harmonicity of one coordinate
function cann't be extended on the others provided $p\ne 2$. This
means that $\R{n+1}$ is equipped with a specified direction $e$. This
kind of asymmetry is typical for the Minkowski spaces $\R{n+1}_{1}$
with time-axis $Oe$. Another example is the tubular minimal surfaces
(see the definition below) which are Euclidean analogous of the
relative strings in nuclear physics.

We should also mention that equation (\ref{plaplace}) is of great
importance in the nonlinear potential theory \cite{Resh} and the
elliptic type PDE's \cite{Lind}, \cite{Hein}. On the other hand, it
is closely linked to quasiregular mappings (see \cite{HKM} for
detailed information).

In the first part of our paper we discuss the basic facts of
$p$-minimal surfaces. In particular, we show that
the Gauss map of a two-dimensional $p$-minimal surface is
$K(p)$-quasiconformal. In section~\ref{ROTAT} we develop a new method
for studying the shape of $p$-minimal tubes. We
establish also the estimates for the sizes of cross-sections of such surfaces
which provides us by information about the evolution of a $p$-minimal tube in the
`time-direction'.

It is also the aim of this paper to demonstrate the properties which are
common for  tubular minimal and $p$-minimal surfaces. We only
briefly discuss the nonparametric case for $p$-minimal surfaces.
The general theory of $p$-minimal surfaces and some examples
will be given in a forthcoming paper.

\section{Preliminary properties of $p$-minimal surfaces}
\label{PRELIM}

{\bf 2.1.}
We have noticed above that the case $p=2$ corresponds to the minimal
surfaces. To clarify the geometrical meaning of (\ref{plaplace})
for arbitrary $p$ we denote by $k_e(m)$ the curvature of $\M$ in
$e$-direction (i.e. the curvature of the section of $\M$ by
2-plane spanning on $e$ and the unit normal $\nu$ to $\M$ at a
point $m$).

\begin{proposition}
Let $m$ be a noncritical point, i.e. $e\wedge \nu(m)\ne 0$. Then
\begin{equation}
H(m) = -(p-2) k_e(m).
\label{geom}
\end{equation}
\label{pr1}
\end{proposition}

\begin{proof}
Really, let $\D$ and $\nabla$ denote the standart
covariant derivatives in $\R{n+1}$ and $\M$ respectively. Then
$$
\nabla f(m) = \left(\D\scal{x(m)}{e}\right)^\top = e^\top,
$$
where $e^\top$ is the projection of $e$ onto the tangent space to the
surface $\M$ at a point $m$. It follows from the assumtions of the
proposition that $|e^\top| \ne 0$ or, it is the same $|\nabla
f(m)| \ne 0$. Thus, for any tangent vector $X$ we obtain
$$
\nabla_X|\nabla f| =\nabla_X|e^\top| =
\frac{\scal{\D_Xe^\top}{e^\top}}{|e^\top|} =
\frac{\scal{\D_X(-e^\bot)}{e^\top}}{|e^\top|} =
\scal{e}{\nu}\scal{A(X)}{\frac{e^\top}{|e^\top|}}.
$$
Here $A$ is the Weingarten map of $\M$ and $e^\bot$ is the projection
of $e$ onto the normal space to $\M$. By virtue of the
symmetry of $A$ we conclude that
\begin{equation} \nabla |\nabla f| =
\scal{e}{\nu}A(\tau),
\label{Wein}
\end{equation}
where $\tau =e^\top/|e^\top|$ is well defined at $m$. After
substituting (\ref{Wein}) into (\ref{plaplace}) we have
$$
\Delta_pf = {\rm
div}(|\nabla f|^{p-2}\nabla f) = (p-2)|\nabla
f|^{p-3}\scal{\nabla f}{\nabla(|\nabla f|)} +
$$
\begin{equation}
+ |\nabla f|^{p-2}\Delta f =
|\nabla f|^{p-4} \scal{e}{\nu}\left ( |e^\top|^2 \Delta f +
(p-2)\scal{A(e^\top)}{e^\top} \right ).
\label{form}
\end{equation}
Finally, the definition of $k_e(m)$ together with the well known
connection between the mean curvature $H(m)$ and the inner Laplacian
\cite{KN}: $\Delta f(m)=H(m)\scal{e}{\nu}$ yield from (\ref{form})
\begin{equation}
\Delta_pf = |\nabla f|^{p-4}\scal{e}{\nu}|e^\top|^2\biggl ( H(m) +
(p-2) k_e(m)\biggr ),
\label{hyper}
\end{equation}
everywhere in the regular part $M_0\equiv\{m\in M:|e^\top(m)|\ne
0\}$.

We assume now that equality (\ref{geom}) doesn't hold at some
noncritical point $m_1 \in M_0$. Then in view of (\ref{hyper}) and
(\ref{plaplace}), $\scal{e}{\nu(m_1)} \equiv 0$, and by continuity of
the expression in parentheses in (\ref{hyper}), the last identity
holds everywhere in some neighbourhood $\Omega(m_1)$. Thus,
in $\Omega(m_1)\cap M_0$ the coordinate function
$f(m)=\scal{e}{x}$ is constant and, it follows that $A\equiv 0$ in
$\Omega(m_1)$. But this conclusion trivially yields validity of
(\ref{geom}) which contradicts our assumption.

\end{proof}

The following assertion is an immediate consequence of the Meusnier
theorem.

\begin{corollary}
The mean curvature $H$ of a $p$-minimal surface $\M$ and the mean
curvature $h$ of the section $\Sigma(\tau)$ linked by
\begin{equation}
h(m) = - \frac{p-1}{\omega} \, k_e(m) = \frac{p-1}{p-2} \;
\frac{H(m)}{\omega}
\label{mene}
\end{equation}
where $\omega=\scal{\nu_m}{e}$.
\label{menje}
\end{corollary}

We use further the auxiliary assertion which clarifies the
local structure of a $p$-minimal surface near a critical point. We
notice that this property has not an analogue in minimal surfaces
theory.

\begin{lemma}
\label{zero}
Let $\M$ be a $p$-minimal surface given as a graph of $C^2$-function
$f(x)$ defined in a domain $G\subset\R{n}$. Let $x_0\in G$ be a
critical point of $f(x)$, i.e. {$\D f(x_0)=0$.} Then the Hessian
$\D^2 f$ is degenerate. Other words, $x_0$ is a planar point.
\end{lemma}

\begin{proof}
To prove this assertion we rewrite (\ref{geom}) in a more suitable
way. In the local coordinates we have the following formulas for the
mean curvature $H(m)$ and the Laplace-Beltrami operator $\Delta$
respectively
$$
H(m) = \frac{1}{g^{3/2}}\sum_{i,j=1}^{n} (g\delta_{ij} - \D_if \D_jf)
\D_{ij}^2f,
$$
\begin{equation}
\Delta u = \frac{1}{\sqrt{g}}\sum_{i,j=1}^{n} \D_j(g^{ij}\sqrt{g}\D_j u),
\label{belt}
\end{equation}
where $\D_i$ denotes the covariant derivative along the coordinate
vector $e_i$, $g^{ij}$ is the inverse matrix to the metric tensor
$g_{ij}=\delta_{ij}+\D_i f\D_j f$ and $g=\det \|g_{ij}\|$. Hence, we
obtain from (\ref{plaplace}) and (\ref{form})
\begin{equation}
g|\D f|^2\;{\rm tr} \D^2 f +
\sum_{l,s=1}^{n} (p-2-|\D f|^2)\D_lf\D_sf \D_{ls}^2f = 0,
\label{eq}
\end{equation}
Here $\D^2f$ is the Hessian of $f(x)$ and the trace ${\rm tr}\D^2f$
is equal to the euclidean Laplace operator in $\R{n}$. Write
$a_{ij}=\D_{ij}^2(x_0)$ and  $A=||a_{ij}||$.
Then for an appropriate choice of $\varepsilon >0$ and every vector
$y\in\R{n}$ such that $|y|<\varepsilon$ we have
$$
\D_kf(x_0+y) = \sum_{i=1}^{n} a_{ki}y_i + o(|y|),
$$
and
$$
|\D f(x_0+y)|^2 = O(|y|^2).
$$

Substituting these relations in (\ref{eq}), we arrive at
$$
\sum_{k,l,s=1}^{n}
\sum_{i,j=1}^{n}
(a_{ki}a_{kj} {\rm tr}A + (p-2)a_{li}a_{sj}a_{ls}) y_iy_j = o(|y|^2).
$$

Taking into account validity of the last equality for all
sufficiently small $y\in\R{n}$ we obtain a matrix equation
\begin{equation}
A^2(I\, {\rm tr}A + (p-2)A ) = 0,
\label{vos}
\end{equation}
where $I$ is the unit matrix. By virtue of symmetry of the Hessian
$A$ we can choose the orthonormal basis of $\R{n}$ consisting of the
eigenvectors of $A$. Namely, $A$ takes a diaganal form
$\lambda_i\delta_{ij}$ and from (\ref{vos}) we have for $i: 1\leq i\leq n$,
$$
\lambda_i\Bigl (\lambda_i(p-2) + {\rm tr}A\Bigr ) = 0.
$$
We see from the last identity that all non-zero eigenvalues
$\lambda_i$ must be equal to {$-(p-2)^{-1}{\rm tr}A$}.
Let $\lambda_1,\ldots,\lambda_k$ be all such numbers. Then after
summing we obtain
\begin{equation}
{\rm tr}A = \sum_{i=1}^{k} \lambda_i = -\frac{k}{p-2}{\rm tr}A.
\label{trace}
\end{equation}
On the other hand ${\rm tr}A=k\lambda_1 \ne 0$. It follows from
(\ref{trace}) that $p=2-k$, where $k\geq 1$ is a positive integer.
But it contradicts with $p > 1$ and hence, all of $\lambda_i $ are
zeroes. Now the theorem follows from the standart properties of
symmetric matrices.
   \end{proof}

\medskip
{\bf 2.2.}
Given a surface $\M$ in $\R{3}$ we denote by $\gamma(m):\M\to S^2$
the standart Gauss map. A result of Gauss state that, if the surface
is minimal that map is conformal. Here we extend this property on
$p$-minimal surfaces. First we remind

\medskip
{\bf Definition 3.} (\cite{Ahlf}, \cite{Resh}) \
A map $F:M_1 \to M_2$ of two smooth Riemannian manifolds
$M_1$ ¨ $M_2$ is called a {\it quasiconformal map}
if the Jacobian $\det d_xF$ doesn't change the sign on $M_1$  and
for almost every $x\in M_1$,
\begin{equation}
\max |d_xF(E)|\leq K_m \min |d_xF(E)|
\label{resh}
\end{equation}
where $\min$ and $\max$ are given over all unit tangent vectors
$E$ of $T_xM_1$. The number $K=\max_{m\in M_1}K_m$ is called the
{\it distortion coefficient} of $F$.

\begin{theorem}
Let $\M$ be a twodimensional $p$-minimal surface in $\R{3}$.
Then the Gauss map is $K(p)$-quasiconformal map with the distortion
coefficient
\begin{equation}
K(p) = \max\{p-1 ; 1/(p-1)\}.
\label{distortion}
\end{equation}
\label{Gauss}
\end{theorem}

\begin{proof}
We notice that the tangent spaces $T_mM$ to $M$ and $T_{\gamma
(m)}S^2$ to the unit sphere $S^2$ can be regarded as canonical
isomorphic ones. Really, we identify the vector $A(E)$ with
$d\gamma_m(E)$, where $d\gamma_m$ is the differential of the Gauss
map at $m$. We specify a point $m\in M$ and choose the orthonormal
basis $E_1, E_2$ of the tangent space $T_mM$ which diagonalizes
$A$, i.e.
$$
A(E_i) = \lambda_i E_i.
$$
Here $\lambda_1,\,\lambda_2$ are the principal curvatures of $\M$
at $m$. Without loss of generality we can arrange that $|e^T(m)|\ne
0$. Really Lemma \ref{zero} yields that the homomorphism $A$ is
the identical zero and (\ref{resh}) is trivial.

Let us denote $\tau = e^T/|e^T|$. Then for some angle $\psi \in [0;\,
2\pi ]$,
$$
\tau = E_1 \cos\psi + E_2\sin\psi ,
$$
and by the Meusnier theorem we have
$$
\scal{A \tau}{\tau} =
\lambda_1\cos^2\psi + \lambda_2\sin^2\psi =
- \frac{1}{p-2}(\lambda_1+\lambda_2).
$$
Hence
$$
\lambda_1 =- \lambda_2 \frac{1+(p-2)\sin^2\psi}{1+(p-2)\cos^2\psi}.
$$

It is a direct consequence of the last identity that the Jacobian
$\det (d_m\gamma )= \lambda_1\lambda_2$ must be negative.
Using standart facts of the quadratic forms theory allowed us to
conclude that the distortion coefficient of $\gamma$ at a point $m$
is less or equal to
$$
K_m=\max_{\psi}\{q ; \frac{1}{q}\} ,
\qquad q=\frac{1+(p-2)\sin^2\psi}{1+(p-2)\cos^2\psi}.
$$
Then varying $\psi$ we obtain the required maximum of $K_m$.
   \end{proof}

L.~Simon in \cite{SL1} established that every entire
twodimensional nonparametric surface with quasiconformal Gauss map
must be plane. As a consequence of this result we obtain a version
of the well-known Bernstein theorem.

\begin{corollary}
\label{cor1}
Let $\M$ be an entire  $p$-minimal graph in $\R{3}$.
Then $\M$ is a plane.
\end{corollary}

{\bf Remark.}
As follows from \cite{MT1} that every minimal
$n$-dimensional graph $\M$ in $\R{n+1}$ has parabolic conformal type.
Other words, every compact on $\M$ have zero $n$-capacity. In
these papers we apply the quasiconformal mapping theory to minimal
surfaces.  The methods used there allows us to conclude that the
similar property holds for $p$-minimal graphs also if $p\geq n$.
These facts together with Corollary \ref{cor1} make very likely to be
true the following:

\medskip
{\bf Conjecture.} Let $\M$ be an entire $p$-minimal graph given over
the whole $\R{n}$. If $p\geq n$ then $\M$ is a hyperplane.

\section{Tubular $p$-minimal hypersurfaces}
\label{ROTAT}

\smallskip
{\bf 3.1.}
In this section we deals with tubular type $p$-minimal
surfaces. This class of surfaces in twodimensional case was involved
by J.~C.~C.~Nitsche \cite{Ni} and have been studied by
V.~M.~Miklyukov \cite{M1} in highdimensional situation.

\medskip
{\bf Definition 4.}
We say that a surface $\M$ is a {\it tube} with the projection
interval $\tau (\M )\subset Ox_{n+1}$, if

(1) for any $\tau\in\tau(\M )$ the
sections $\Sigma_\tau = x(M)\cap \Pi_\tau$ by hyperplanes
$\Pi_\tau =\{ x\in \R{n+1}_1 : x_{n+1}=\tau\}$ are not empty
compact sets;

(2) for $\tau', \tau'' \in \tau(\M)$ any part of $\M$ situated
between two different $\Pi_{\tau'}$ and $\Pi_{\tau''}$ is a
compact set.

If $\tau(\M)$ is an infinite interval we call the surface to be an
{\it infinite tube}. Otherwise, we call a length of $\tau(\M)$ the
{\it life-time} of $\M$.

Let
$$
\rho(\tau)=\max_{m\in
\Sigma(\tau)}\left(\sum_{i=1}^{n}x_i^2\right)^{1/2}.
$$
It follows from the results of \cite{VM}, \cite{MT2}, \cite{Kl}, that
 every $n$-dimensional minimal tube of arbitrary
codimension satisfies the following differential inequality
\begin{equation}
\rho(\tau)\rho''(\tau)\geq (n-1)(1+\rho'(\tau)^2),
\label{veden}
\end{equation}
which is crucial for all theory of minimal tubes. As a consequence
every minimal tube for $n\geq 3$ contains in a slab between two
parallel planes. Hence, there are not many-dimensional infinite
minimal tubes. In contrast, the twodimensional case essentially
differs from the highdimensional one: there are tubes of finite
life-time as well as infinite tubes. Moreover, we shown in \cite{T3}
that the life-time in the first case is derived by the angle between
the full-flow vector of a minimal tube and the time-axe.

\begin{lemma}
Let $V$ be a convex compact in $\R{n}$ and $W$ be a compact such that
$W\setminus V\ne \emptyset$. Then there exists a closed ball $B\subset
\R{n}$ such that
\begin{equation}
W\subset B
\label{eq21}
\end{equation}
and
\begin{equation}
\partial B\cap (W\setminus V)\ne \emptyset.
\label{eq22}
\end{equation}
\label{ball}
\end{lemma}

\begin{proof}
The distance function $f(x)={\rm dist}(x, V)$ is continuous on
$\R{n}$.  It follows from the conditions of the lemma that this
function attains the maximum value on $W$ at some point $a\in W$ and
$d=f(a)>0$. On the other hand, by virtue of convexity
of $V$ there exist the unique point $b\in \partial V$ such that
$f(a)=\|b-a\|$.

We choose a new coordinate system of $\R{n}$ with the origin at
$a$, the first coordinate vector
$$
e_1=\frac{b-a}{d}
$$
and the others $e_2,\ldots,e_n$ to be based an orthonormal system
together with $e_1$. Then, the hyperplane given by $x_1=d$ is one of
support to $V$ at $a$. It follows from the triangle inequality that
$W$ contains in halfspace $\{x_1\geq 0\}$ and $V$ in $\{x_1\geq d\}$.

Given positive $h$ and $R$ we specify an open ball
$$
B(R,h)=\{x\in\R{n}:(x_1+R)^2+x_2^2+\ldots x_n^2<(R+h)^2\}.
$$
By our choice and compactness of $V$, given a positive $\epsilon$
there exists $R>0$ such that $V$ containes in the ball $B(R,\epsilon)$.

Suppose $\epsilon=d/2$ and $R_0$ is the corresponding radius.
Then the definition of $d$ yields that $a\not\in B(R_0,d/2)$, however
the greater ball $B(R_0,3d/2)$ contains $V$ as well as $W$. Let
$\delta_0$ is the minimum over all $\delta\in (0;d)$ such that
$$
W\subset \overline{B(R_0,d/2+\delta)}.
$$
Then $a\in\partial B$, where $B=\overline{B(R_0,d/2+\delta_0)}$ and
$V\cap B=\emptyset$.

   \end{proof}

\begin{corollary}
[Maximum Principle]
Let $\M=(M,x)$ be an immersed compact $p$-minimal hypersurface in
$\R{n+1}$ with nonempty boundary $\partial M$. Then
\begin{equation}
\conv x(\partial M)=\conv x(M),
\label{max}
\end{equation}
where $\conv E$ is the convex hull of $E$.
\label{maximum}
\end{corollary}

\begin{proof}
Let us denote $\Omega=\conv x(\partial M)$ and assume that
(\ref{max}) fails. Then it implies $x(M)\setminus \Omega\ne
\emptyset$. By Lemma \ref{ball} we can find the closed ball $B$ such
that $x(M)\subset B$ and there exists a point $m\in {\rm int}\; M$,
$x(m)\in\partial B$. We choose the neighbourhood $\mathcal{O}$ of $m$
such that the restriction of $x$ on $\mathcal{O}$ is embedding. The
further arguments will be local and we can arrange that $\M=x(\mathcal{O})$ without loss of generality.

Because of the choice of $B$, the tangent spaces to $\M$
and $\partial B$ at $x(m)$ coincide. Moreover, $\M\subset B$ and the
standart comparison principle for touching surfaces gives the
following inequality
\begin{equation}
\lambda_i\geq\frac{1}{R},
\label{comparison}
\end{equation}
where $\lambda_i$ are the principal curvatures of $\M$ at $m$ with
respect to the inward normal of $\partial B$ and $R$ is the radius of
$B$.

We now turn to identity (\ref{geom}). By the definition of $k_e(m)$
there exists a system of positive numbers $\alpha_i\leq 1$ such that
$$
\sum_{i=1}^{n}{\alpha_i}=1
$$
and
$$
k_e(m)=\sum_{i=1}^{n}\alpha_i\lambda_i
$$
It follows from these relations, (\ref{comparison}) and (\ref{geom})
that
$$
0=\sum_{i=1}^{n}\lambda_i(1+(p-2)\alpha_i)\geq \frac{n+p-2}{R}>0.
$$
This contradiction proves Corollary \ref{maximum}.
   \end{proof}

\smallskip
{\bf 3.2.}
Further we use the Minkowski operations. Namely, given $A,B\subset
\R{n}$ the notations $A\oplus B$ and $\lambda A$ reserved for the
sets $\{x=a+b: a\in A, b\in B\}$ and $\{x=\lambda a: a\in A\}$.

{\bf Definition 5.}
A family of convex sets $\{\Omega(\tau):\tau\in[\alpha,\beta]\}$ is
called \cite{Leicht} {\it convex} if for arbitrary $\tau_1<\tau_2$
from interval  $[\alpha;\beta]$ and a nonnegative $t\leq 1$ one holds
$$
\Omega\left(\tau_1 t+\tau_2 (1-t)\right)\subset
t\Omega(\tau_1)\oplus \overt\Omega(\tau_2).
$$

Let $\M$ be an $n$-dimensional $p$-minimal tube in $\R{n+1}$. Let
us denote by $\Omega(\tau)$ the projection of the convex hull of the
section $\Sigma(\tau)$ onto the hyperplane $\Pi_0=\{x_{n+1}=0\}$.
Then
$$
\conv \Sigma(\tau)=\tau e_{n+1}\oplus\Omega(\tau).
$$

\begin{theorem}
The family $\{\Omega(\tau):\tau\in \tau(\M)\}$ is convex.
\label{convex}
\end{theorem}

\begin{proof}
We specify $\tau_1<\tau_2$ from interval $\tau(\M)$ and $t\in
[0;1]$. Let $H$ be the slab $\{x: x_{n+1}\in(\tau_1;\tau_2)\}$ and
$M'=x^{-1}(H\cap x(M))$.
Then Corollary \ref{maximum} gives
$$
V\equiv \conv (\Sigma(\tau_1)\cup \Sigma(\tau_2))= \conv x(M').
$$

Let $\tau_0=t\tau_1+\overt \tau_2$. Then $\Sigma(\tau_0)\subset V$
and by the definition of the convex hull we conclude that
$\conv\Sigma(\tau_0)\subset V$.

We choose an arbitrary $z\in \Omega(\tau_0)$. Then $y=z+\tau_0
e_{n+1}\in \Pi(\tau_0)\cap V$, and there exist
$y_i\in\conv\Sigma(\tau_i)$ and $\lambda\in [0;1]$ such that
\begin{equation}
y=\lambda y_1+\overline{\lambda} y_2.
\label{igrek}
\end{equation}
Then decomposition of $y_i=z_i+\tau_i e_{n+1}$ for certain
$z_i\in\Omega(\tau_i)$ and (\ref{igrek}) give
$$
z=\lambda z_1+\overline{\lambda} z_2, \qquad \tau_0=
\lambda \tau_1+\overline{\lambda} \tau_2.
$$
Hence, $\lambda=t$ and it follows $z\in t\Omega(\tau_1)\oplus
\overt\Omega(\tau_2)$ as required.
   \end{proof}

The following assertion gives a sample of applications of the last
result.

\begin{corollary}
Let $R(\tau)$ be the radius of the least ball which contains
$\Sigma(\tau)$ {\rm (further we call such a ball to be circumscribed
near $\Sigma(\tau)$)}. Then $R(\tau)$ is a convex function.
\label{radius}
\end{corollary}

\begin{proof}
We denote by $B(\tau)$ the projection of the circumscribed
near $\Sigma(\tau)$ ball onto $\Pi_0$. Then, by virtue of convexity
of $B(\tau)$ we have $B(\tau)\supset\Omega(\tau)$, and Theorem
\ref{convex} yields for arbitrary $t\in[0;1]$
$$
\Omega(\tau_0)\subset
t\Omega(\tau_1)\oplus \overt \Omega(\tau_2)\subset
t B(\tau_1)\oplus \overt B(\tau_2)=B_0,
$$
where $\tau_0=\tau_1 t+\tau_2 (1-t)$. By the definition, $R(\tau_0)\leq
R_0$, where $R_0$ is the radius of $B_0$. On the other hand, $R_0=t
R(\tau_1)+ \overt R(\tau_2)$ and we obtain the required inequality
$$
R(\tau_1 t+\tau_2 (1-t))\leq R(\tau_1)+ \overt R(\tau_2).
$$
   \end{proof}

\smallskip
{\bf 3.3.}
Now we study the structure of $\Sigma(\tau)$ in more detailes. This
requires further delicate information not only about $R(\tau)$ but
about the curve of the centers of the balls $B(\tau)$ as well. Let us
denote by $\xi(\tau)$ the center of $B(\tau)$. We remind without proof
the well known extremal property of $B(\tau)$ (see \cite{Leicht},
Theorem~7.5).

\begin{lemma}
Let $E$ be a closed subset of $\R{n}$ and $B(E)$ the circumscribed
near $E$ ball with the center $\xi$. Then for all unit vectors
$y\in\R{n}$ there exists $b\in \partial E\cap \partial B(E)$ such
that
\begin{equation}
\scal{b-\xi}{y}\geq 0.
\label{scal}
\end{equation}
\label{extrem}
\end{lemma}

We denote by
$$
\sigma(E)=\min_{y\in S^{n-1}}\max_{b\in \partial B\cap
E}\frac{\scal{b-\xi}{y}}{R},
$$
where $B$ is the circumscribed ball near a compact $E$, $R$ is the
radius and $\xi$ is the center of $B$. It follows from (\ref{scal})
that $0\leq \sigma(E)\leq 1$. Moreover, one easy to see that
$\sigma(E)=0$ if and only if the intersection of the boundary sphere
$S=\partial B$ with $F$ lies in some equatorial semisphere of $S$.

\begin{theorem}
Let $\M$ be a $p$-minimal tube in $\R{n+1}$ such that
\begin{equation}
\sigma(\Sigma(\tau))\geq \epsilon>0, \qquad \forall \tau\in\tau(\M).
\label{condition}
\end{equation}
Then $\xi(\tau)$ is a $\delta$-convex curve of $\tau$. Other words,
any coordinate function $\xi_k(\tau)$ admits the composition
$$
\xi_k(\tau)=\varphi_k(\tau)-\psi_k(\tau),
$$
with $\varphi_k(\tau)$, $\psi_k(\tau)$ to be convex functions.
\label{sigma}
\end{theorem}

\begin{proof}
We consider $\tau_1$, $\tau_2$ from $\tau(\M)$ and $t\in[0;1]$. Let
us denote by $B(\tau_i)=B_i(\xi(\tau_i),R_i)$ the corresponding
circumscribed near $\Sigma(\tau_i)$ balls. As above we have for
$\tau_0=t \tau_1+\overt \tau_2$
$$
\Omega(\tau_0)\subset
t B(\tau_1)\oplus \overt B(\tau_2).
$$
In force of Lemma \ref{extrem} we can find $y\in \partial B(\tau_0)\cap
\Sigma(\tau_0)$ such that
$$
\scal{y-\xi(\tau_0)}{\xi(\tau_0)-\xi_0}\geq
\epsilon |y-\xi(\tau_0)|\cdot |\xi(\tau_0)-\xi_0|,
$$
where $\xi_0=t \xi(\tau_1)+\overt \xi(\tau_2)$. Hence,
$$
|y-\xi_0|^2=\bigl|(y-\xi(\tau_0))+(\xi(\tau_0)-\xi_0)\bigr|^2\geq
$$
$$
\geq |y-\xi(\tau_0)|^2+|\xi(\tau_0)-\xi_0|^2+2\epsilon
|y-\xi(\tau_0)|\cdot|\xi(\tau_0)-\xi_0|,
$$
and taking into account that $|y-\xi(\tau_0)|=R(\tau_0)$ and
$|y-\xi_0|\leq R_0$ we obtain
$$
|\xi(\tau_0)-\xi_0|^2+2\epsilon |y-\xi(\tau_0)|\cdot|\xi(\tau_0)-\xi_0|+
(R^2(\tau_0)-R^2_0)\leq 0,
$$
and as a consequence,
\begin{equation}
|\xi(\tau_0)-\xi_0|\leq
\frac{R^2_0-R^2(\tau_0)}{R(\tau_0)\epsilon
+\sqrt{R^2_0-R^2(\tau_0)(1-\epsilon^2)}}.
\label{root}
\end{equation}
By Corollary \ref{radius} we have $R_0\geq R(\tau_0)$ and from
(\ref{root}),
\begin{equation}
|\xi(\tau_0)-\xi_0|\leq
\frac{R^2_0-R^2(\tau_0)}{\epsilon(R(\tau_0)+R_0)}
=\frac{1}{\epsilon}(R_0-R(\tau_0)).
\label{delta}
\end{equation}

We consider the coordinate function
$\xi_k(\tau)=\scal{\xi(\tau)}{e_k}$. Then (\ref{delta}) yields
$$
t\xi_k(\tau_1)+\overt\xi_k(\tau_2)-\xi_k(\tau_0)\leq
\frac{1}{\epsilon}\bigl(
tR(\tau_1)+\overt R(\tau_2)-R(\tau_0)
\bigr).
$$
This inequality means that the difference $\psi(\tau)=\epsilon^{-1}
R(\tau)-\xi_k(\tau)$ is convex. Therefore, by Corollary \ref{radius}
we obtain the required decomposition of $\xi_k(\tau)$ into
difference of two convex functions
$$
\xi_k(\tau)=\frac{1}{\epsilon}R(\tau)-\psi(\tau)
$$
and the lemma is proved.
   \end{proof}

\begin{theorem}
Let $\M$ be a $p$-minimal surface with assumption {\rm
(\ref{condition})} and $\beta=(n-1)/(p-1)$. Then $R(\tau)$ and
$\xi(\tau)$ satisfy the differential inequality
\begin{equation}
R(\tau)R''(\tau)\geq
\beta(1+R'(\tau)^2)+|\xi'(\tau)|^2 \min\{\beta;1\}
\label{vedennew}
\end{equation}
almost everywhere in $\tau(\M)$.
\end{theorem}

\begin{proof}
Convexity of a function provides existence a.e. of the
second differential (see \cite{Leicht} or \cite{Resh}, Theorem 5.3).
It follows from Corollary \ref{radius}, Theorem \ref{sigma}
that $R(\tau)$ as well as $\xi_k(\tau)$ have the second differentials
almost everywhere in $\tau(\M)$. We denote by $\tau'(\M)$ the set of
full measure where the second differentials of $R(\tau)$ and
$\xi_k(\tau)$, $1\leq k\leq n+1$ do exist.

Let $S^{n-1}$ be the unit sphere in $\Pi_0\sim \R{n}$ endowed by the
standart metric. We consider the hypersurface $\M_0$ given by
$$
w(\theta,\tau)=\xi(\tau)+R(\tau)\theta+\tau e_{n+1}\;
:S^{n-1}\times \R{}\rightarrow \R{n+1}
$$
where $\theta\in S^{n-1}$.
We have shown in \cite{T2} that for such a surface the curvature
$k_{e,\M_0}$ in $e$-direction can be calculated by
\begin{equation}
k_{e,\M_0}(\theta,\tau)=
\frac{\omega^3}{R(\tau)}\left[
R(\tau)R''(\tau)+R(\tau)\scal{\xi''(\tau)}{\theta}+
\scal{\xi'(\tau)}{\theta}^2-|\xi'|^2
\right]
\label{k}
\end{equation}
where
$$
\omega^2=\scal{\nu_m}{e}^2=
\frac{1}{1+\biggl(R'(\tau)+\scal{\theta}{\xi'(\tau)}\biggr)^2}.
$$

By the definition of functions $R(\tau)$ and $\xi(\tau)$ we conclude
that the surface $\M$ contains inside of $\M_0$ in the sense that
$\Sigma(\tau)$ is a subset of $\Pi(\tau)\cap \M_0$ for all
$\tau\in\tau(\M)$.

Let us consider an arbitrary $\tau\in\tau'(\M)$ and $E=\Omega(\tau)\cap
\partial B(\tau)$. The surfaces $\M$ and $\M_0$ have the common
outward normal $\nu_m$ at $m=y\oplus \tau e_{n+1}$ for every $y\in
E$ (we mean by {\it outward} the normal which is directed out from
the inside of $B(\tau)$). Let $\mathcal{O}$ be the neighbourhood of $m$
where $x(\cdot)$ is an embedding. It is a consequence of the
definition of $\M_0$ that $\nu_m\wedge e_{n+1}\ne 0$. We denote by
$\gamma(\tau)$ and $\gamma_0(\tau)$ the sections of $x(M)$ and $\M_0$
by the two-plane spanning on $\nu_m$ and $e_{n+1}$. Then the
comparison principle for touching surfaces yields
$$
k_{e,\M}(m)\leq k_{e,\M_0}(m).
$$

We write $h(m)$ and $h_{0}(m)$ for the mean curvatures at $m$
of the sections $\Sigma(\tau)$ and $\Pi(\tau)\cap\M_0=\xi(\tau)\oplus
\tau e_{n+1}\oplus B(\tau)$ with respect to their common outward
normal. Then the comparison principle arrive at the inequality
$$
h(m)\leq h_{0}(m)\equiv -\frac{n-1}{R(\tau)},
$$
and after (\ref{mene})
$$
- \frac{p-1}{\omega} \, k_e(m) \leq
 -\frac{n-1}{R(\tau)}.
$$
By (\ref{k}) we obtain after simplification
\begin{equation}
R(\tau)R''(\tau)-\beta(1+R'(\tau)^2)\geq
(\beta-1)\scal{\xi'(\tau)}{\theta}^2
+|\xi'|^2+\scal{\theta}{y},
\label{q}
\end{equation}
where $y=2\beta R'(\tau)\xi'(\tau)-R(\tau\xi''(\tau))$.
Thus, Lemma \ref{extrem} to be applied to the vector $y$ provides
$b\in E$ such that $\scal{b-\xi(\tau)}{y}\geq 0$. We take
$$
\theta_0=\frac{b-\xi(\tau)}{R(\tau)}
$$
and it follows from (\ref{q})
$$
R(\tau)R''(\tau)-\beta(1+R'(\tau)^2)\geq
(\beta-1)\scal{\xi'(\tau)}{\theta_0}^2
+|\xi'(\tau)|^2\geq |\xi'(\tau)|^2 \min\{\beta;1\},
$$
and the theorem is proved completely.

   \end{proof}

{\bf Remark.}
Finally, we notice that the quantity $R(\tau)$ measures the size of
the section $\Sigma(\tau)$ instead of the distance this section from a
fixed line in the previous inequalities (\ref{veden}). Moreover, in
the base case $p=2$ the established inequality (\ref{vedennew}) is
more strong then (\ref{veden}).

On the other hand, $\delta$-convex functions belong to the class
$\overline{W}^2_{1,{\rm loc}}(\tau(\M))$; that is, has a second-order
generalized derivative that is a measure (see \cite{Resh}, Chapter 2,
\S 4.10, Corollary). This allows to proceed the integration of
(\ref{vedennew}) to comletion in the standart way \cite{VM}:

\begin{corollary}
Let $\M$ be a $p$-minimal tube, ${\rm dim}\M=n>p>1$. Then has finite
life-time ${\rm length}\tau(\M)$. Moreover,
$$
{\rm length}\tau(\M)\leq 2c_\beta r(\M), \quad \beta=\frac{n-1}{p-1}
$$
where
$$
r(\M)\equiv \min_{\tau\in\tau(\M)}R(\tau)>0
$$
and
$$
c_{\beta}=\int_{0}^{+\infty}\frac{d\,t}{(1+t^{2\beta})^{1/2}}.
$$
\end{corollary}

\begin{small}

\end{small}


\begin{thebibliography}{99}
\bibitem{Ahlf}
L. Ahlfors,
{\it Lectures on quasiconformal mappings},\/ Toronto-New York-London:
Van Nostrand Math. Studies, 1966.
\bibitem{Resh}
V. M. Goldstein and Yu.G.Reshetnyak, {\it Introduction in the theory of
function with generelized derivatives and quasyconformal mappings},\/
Nauka, Moscow, 1983.
\bibitem{HKM}
J. Heinonen, T. Kilpelainen and  O. Martio,
{\it Nonlinear potential theory of degenerate elliptic equations},\/
Oxford Univ. Press, London, 1993.
\bibitem{Hein}
J. Heinonen, T. Kilpelainen and  J. Rossi, {\it The grows of
$A$-subharmonic functions and quasiregular mappings along asymptotic
paths},\/ Indiana Univ. Math. J., {\bf 38}(1989), N3,  581-601.

\bibitem{Kl}
V. A. Klyachin,
{\it Estimate of spread for minimal surfaces of arbitrary
codimension}\/, Sibirsk. Mat. Zh.,  {\bf 33}(1992), N 5, p.201-207.

\bibitem{KN}
S. Kobayashi and K. Nomizu, {\it Foundations of differential
geometry},\/ Vol.2, Interscience, 1969.

\bibitem{Lind}
P. Lindqvist, {\it On the definition and properties of
$p$-superharmonic functions},\/ J. Reine Angew. Math., {\bf
365}(1986), 67-75.

\bibitem{Leicht}
K. Leichtweiss, {\it Konvexe Mengen}\/, Springer, 1980.

\bibitem{M1}
V. M. Miklyukov,
{\it On some properties of tubular minimal surfaces in $R^n$}\/,
Dokl. Akad. Nauk SSSR, {\bf 247}(1979), No.3, 549-552; English
transl. in Soviet Math. Dokl. {\bf 20}(1979).

\bibitem{MT1}
V. M. Miklyukov and V. G. Tkachev, {\it On the structure in the large
of externally complete minimal surfaces in $\R{3}$}\/, Sov. Math.
Izv. VUZ, {\bf 31}(1987), 30-36.

\bibitem{MT2}
V. M. Miklyukov and V. G. Tkachev, {\it Some properties of tubular
minimal surfaces of arbitrary codimension}\/, Mat. Sb. {\bf
180}(1989), N 9, 1278-1295; English transl. in Math. USSR Sb. {\bf
68}(1991) No. 1, 133-150.

\bibitem{VM}
V. M. Miklyukov and V. D. Vedenyapin,
{\it Extrinic dimension of tubular minimal hypersurfaces}\/,
Mat. Sb., {\bf 131}(1986), 240-250; English transl. in Math. USSR Sb.
{\bf 59}(1988).

\bibitem{Ni}
J. C. C. Nitsche, {\it A uniqueness theorem of Bernstein's type for
minimal surfaces in cylindrical coordinates}/\, J. of Math. and Mech.
{\bf 11}3(1962), 293-302.

\bibitem{Os}
R. Osserman, {\it A survey of minimal surfaces},\/ New York: Dover
Publications, 1987.


\bibitem{SL1}
L. Simon, {\it A Holder estimate for quasiconformal mappings between
surfaces in Euclidean space, with application to graphs, having
quasiconformal Gauss map},\/ Acta math. {\bf 139}(1977), 19-51.

\bibitem{SL2}
L. Simon, {\it Equation of mean curvature type in two independed
variables},\/ Pacif. J. Math.  {\bf 69}(1977), N1, 245-268.

\bibitem{T2}
V. G. Tkachev, {\it The external estimates of clasp function for
elliptic hypersurfaces}\/, Preprint 2031-B92, deposited at VINITI,
1992, 17 P. (in Russian)

\bibitem{T3}
V. G. Tkachev, {\it Minimal tubes and coefficients of holomorphic
functions in annulus},\/ in Generalizations of Compl.
Anal. and Appl. in Physics, Centre Nat. de la Recherche Sci., Paris,
1995 {\it to appear}.

\end{thebibliography}
\end{document}